\documentclass{amsart}
\usepackage[T1]{fontenc}

\usepackage{amsfonts, amssymb, amsmath, amsthm, enumerate}
\usepackage{multicol, xspace, array, multirow,tikz-cd}
\usepackage[margin=1in]{geometry}
\usepackage{tcolorbox}
\usepackage[utf8]{inputenc}
\usepackage{multicol}

\theoremstyle{plain}
\newtheorem{theorem}{Theorem}[section]

\theoremstyle{definition}
\newtheorem{definition}[theorem]{Definition}

\newtheorem{lemma}[theorem]{Lemma}
\newtheorem{claim}{Claim}[section]

\title{Existence of a Model of \(o(\kappa)=\kappa^{++}\) from Failure of GCH at a Measurable Cardinal}
\author{Connor Watson}
\date{}

\begin{document}
	\maketitle
	
	\begin{abstract}
		It is well-known that the consistency strength of the GCH failing at a measurable cardinal is the existence of a cardinal \(\kappa\) with \(o(\kappa)=\kappa^{++}\). As the literature does not contain more than a proof sketch of the lower bound of this equiconsistency, we give an expository proof which fills in the details in order to fill this gap in the literature.
	\end{abstract}
	
	\section{Introduction}
	It is known by work of Gitik \cite{GitikUpperBound} and Mitchell \cite{MitchellOGPapers} that the exact consistency strength of the GCH failing at a measurable cardinal is the existence of a cardinal \(\kappa\) of Mitchell order \(\kappa^{++}\). However, the only places where the lower bound of this equiconsistency result is written down are identical proof sketches in \cite{JechSetTheory} and \cite{MitchellCoveringLemma}, and a very outdated proof in \cite{MitchellOGPapers}. As we find that none of the current proofs are approachable for someone who doesn't know the material already, we provide a proof which is essentially the proof given in \cite{JechSetTheory} and \cite{MitchellCoveringLemma} with all of the details filled in. This paper is part of the author's forthcoming Master's thesis, which, among other things, aims to give the first approachable account of Mitchell's core model \(K\) up to \(o(\kappa)=\kappa^{++}\). None of the results presented in this document are originally due to the author.
	
	The author would like to thank Sean Cody for his invaluable comments during the writing of this document.
	
	\section{The Proof}
	We start by establishing conventions. First, we recall the definition of an iterated ultrapower. 
	
	\begin{definition}
		An \textit{iterated ultrapower} of a model \(M\) is any member of a sequence \(\langle M_\alpha:\alpha<\beta\rangle\) constructed inductively on \(\alpha\) as follows:
		
		\(M_0=M\),
		
		\(M_{\alpha+1}=\mathrm{Ult}(M_\alpha,U_\alpha)\) where \(U_\alpha\in M_\alpha\) is a \(\kappa^{(\alpha)}\)-complete nonprincipal ultrafilter on \(\kappa^{(\alpha)}\) for some cardinal \(\kappa^{(\alpha)}\in M_\alpha\),
		
		If \(\gamma\) is a limit ordinal, then \(M_{\gamma}=\varinjlim_{\alpha<\gamma}\langle(M_{\alpha},i_{\alpha,\beta}):\alpha<\beta<\gamma\rangle\), where the maps \(i_{\alpha,\beta}\) are defined in the usual way.
	\end{definition}
	
The reader should note that in forming 	\(M_{\alpha+1}\), we take the ultrapower of \(M_\alpha\) by an ultrafilter which is a member of \(M_\alpha\). There are more general iterated ultrapowers in which ultrapowers are taken by arbitrary ultrafilters in \(V\), but we do not consider these.
	
If \(M_\alpha\) is well-founded, then we identify \(M_\alpha\) with its transitive collapse. Furthermore, we also assume that all iterated ultrapowers are \textit{normal}, i.e. have nondecreasing critical points.\\

We now state the most important properties of \(K\), the core model up to \(o(\kappa)=\kappa^{++}\), which we will need for the lower bound of the equiconsistency we are after. The reader does not need to know the particular about how \(K\) is constructed, just that intuitively, \(K\) should contain all of the ``large cardinal information\("\) below \(o(\kappa)=\kappa^{++}\) which exists in \(V\). A proof sketch of the following lemma can be found in the author's forthcoming master's thesis, and a full proof can be found in \cite{MitchellCoveringLemma}.
		
	\begin{lemma}
		\label{beginninglemma}
		
		\begin{enumerate}
			
		\item If there is no model with a measurable cardinal \(\kappa\) such that \(o(\kappa)=\kappa^{++}\), then every embedding \(i:K\to N\) into a well-founded model \(N\) is an iterated ultrapower by measures in \(K\).
		
		\item If \(U\) is a normal \(K\)-ultrafilter on \(\kappa\) and \(\mathrm{Ult}(K,U)\) is well-founded then \(U\in K\). If \(\mathrm{crit}(U)>\omega_2\) then the hypothesis that \(\mathrm{Ult}(K,U)\) is well-founded is unnecessary.
		
		\end{enumerate}
	\end{lemma}

We are now ready to state and prove the lower bound of the equiconsistency, i.e. that the consistency of the GCH failing at a measurable cardinal \(\kappa\) implies the consistency of the existence of a measurable cardinal \(\lambda\) with \(o(\lambda)=\lambda^{++}\).

	\begin{theorem}
		If there is a measurable cardinal \(\kappa\) with \(2^\kappa>\kappa^+\), then there exists an inner model with a measurable cardinal \(\lambda\) such that \(o(\lambda)=\lambda^{++}\).
		
		\begin{proof}
			Towards a contradiction, suppose \(\kappa\) is measurable and \(2^\kappa>\kappa^{+}\), but there is no inner model with a measurable cardinal \(\lambda\) such that \(o(\lambda)=\lambda^{++}\). In particular, there is no measurable cardinal \(\lambda\) such that \(o(\lambda)=\lambda^{++}\) in \(K\). Even more particularly, \(o(\kappa)\neq\kappa^{++}\) in \(K\). Since we always have that \(o(\kappa)\leq(2^\kappa)^+\), and the GCH holds in \(K\), we have that \(o(\kappa)\leq\kappa^{++}\) in \(K\), so in fact \(o(\kappa)<\kappa^{++}\) in \(K\).
			
			Let \(U\) be any measure on \(\kappa\) and let \(i_U:V\to M=\mathrm{Ult}(V,U)\) be the usual ultrapower embedding. First, we claim that \(|i_U(\kappa)|=2^\kappa\). We have that \(|i_U(\kappa)|\leq 2^\kappa\) since any ordinal \(\alpha<i_U(\kappa)\) is represented by a function \(f:\kappa\to\kappa\), of which there are \(2^\kappa\). To see that \(i_U(\kappa)\geq2^\kappa\), notice that \(\kappa\) is measurable in \(V\), so \(i_U(\kappa)\) is measurable in \(M\). In particular, \(i_U(\kappa)\) is a strong limit cardinal in \(M\). Let \(\mu\) be any cardinal less than \((2^\kappa)^V\). Then \(\mu\) is also less than \((2^\kappa)^M\). Hence any \(V\)-cardinal less than \(2^\kappa\) is not strong limit in \(M\), so \(i_U(\kappa)\geq 2^\kappa\).
			
			Now, let \(i=i_U\upharpoonright K:K\to K^M\) be the restriction of the embedding \(i_U\) to \(K\). By Lemma \ref{beginninglemma}(i), \(i\) is an iterated ultrapower of \(K\), so let \(\langle N_\nu:\nu\leq\theta\rangle\) be the iterates, so that \(N_0=K\) and \(N_{\theta}=K^M\). If \(\nu<\theta\) is a limit ordinal then there are \(\xi_\nu<\nu\) and \(U_\nu\in N_{\xi_{\nu}}\) such that \(N_{\nu+1}=\mathrm{Ult}(N_\nu,i_{\xi_{\nu},\nu}(U_\nu))\). (The reader should be careful and notice that \(U_\nu\notin N_\nu\), but \(i_{\xi_\nu,\nu}(U_\nu)\in N_\nu\), so that \(N_{\nu+1}\) is actually the ultrapower of \(N_\nu\) by a measure in \(N_\nu\).) Moreover, for \(\nu<\theta\) we write \(\kappa_\nu\) for the ordinal \(i_{0,\nu}(\kappa)\). We first claim that this iteration is of length at least \(\kappa^{++}\).
			
			\begin{claim}
				\(\theta\geq\kappa^{++}\).
				\begin{proof}
					We prove the claim by bounding the size of \(i_{0,\nu}(\kappa)\) for \(\nu<\kappa^{++}\) by considering how \(N_\nu\) is represented in terms of the extender of the iteration at that point. In particular,
					
					$$i_{0,\nu}(\kappa)=\{[a,f]:a\in\{\kappa_\eta:\eta<\nu\}^{<\omega}\text{ and }f:[\kappa]^{<\omega}\to\kappa,\;f\in K\}.$$
					
					By simple cardinal arithmetic, \(|\{\kappa_\eta:\eta<\nu\}^{<\omega}|=|\nu|<\kappa^{++}\) and \(\kappa^{(\kappa^{<\omega})}=\kappa^{+}\), so  \(|i_{0,\nu}(\kappa)|\leq|\nu|\cdot\kappa^+<\kappa^{++}\). But we showed previously that \(|i(\kappa)|=|i_{0,\theta}(\kappa)|=2^\kappa>\kappa^{+}\), a contradiction. 
				\end{proof}
			\end{claim}
			
			\begin{claim}
				There is a stationary class \(\Gamma\subseteq\kappa^{++}\) of ordinals of cofinality \(\omega\) such that \(\xi_{\nu}=\bar{\xi}\) and \(U_\nu=\bar{U}\) are constant for \(\nu\in\Gamma.\)
				
				\begin{proof}
					The function \(\nu\mapsto\xi_\nu\) is regressive on limit points of \(\kappa^{++}\) of cofinality \(\omega\), so by Fodor's Lemma there is some \(\Gamma'\subset\mathrm{Lim}(\kappa^{++})\) of points of cofinality \(\omega\) on which \(\nu\mapsto\xi_\nu\) is constant. Denote this constant value \(\bar{\xi}.\)
					
					Let $$\Phi:\Gamma'\to N_{\bar{\xi}}$$ be the function which sends an ordinal \(\nu\) to the measure \(i_{\bar{\xi},\nu}(U_\nu)\). We know that every measure in the iteration has critical point at most \(i(\kappa)\), so \(\mathrm{crit}(\Phi(\nu))\leq\kappa_{\bar{\xi}}\) for every \(\nu\in\Gamma'.\) We now wish to find a bound for the size of \(\mathrm{ran}(\Phi)\), and we do this by counting the number of measures on cardinals in \(N_{\bar{\xi}}\) less than or equal to \(\kappa_{\bar{\xi}}\). Since \(N_{\bar{\xi}}\models |o(\kappa_{\bar{\xi}})|\leq\kappa_{\bar{\xi}}^+\) by our contradiction assumption, and \((\kappa_{\bar{\xi}}^+)^{N_{\bar{\xi}}}<\kappa^{++}\), in \(V\) we can see that \(|\mathrm{ran}(\Phi)|\leq\kappa^+\). Let \(\langle a_\beta:\beta<\alpha\rangle\) where \(\alpha\leq\kappa^+\) be a well-ordering of \(\mathrm{ran}(\Phi)\). Now let $$\Psi:\Gamma'\to\alpha$$ be the function so that \(\Psi(\nu)= a_{\Phi(\nu)}\). Once again, we may use Fodor's lemma to find a stationary \(\Gamma\subseteq\Gamma'\) on which \(\Psi\) is constant. Denote this constant value \(\bar{a}\), and let \(\eta\) be the ordinal so that \(\bar{a}=a_{\Phi(\eta)}\). Then the measure \(\bar{U}=\Phi(\eta)\), the stationary class \(\Gamma\), and the ordinal \(\bar{\xi}\) are as in the statement of Claim 2.2. This completes the proof of the Claim.
				\end{proof}
			\end{claim}
			
			Now fix a particular \(\nu\in\Gamma\cap\mathrm{Lim}(\Gamma)\). Let \(\vec{\kappa}=\langle\nu_n:n<\omega\rangle\) be a cofinal sequence in \(\Gamma\cap\nu\), and for each \(n<\omega\) let \(\kappa^{(n)}=\mathrm{crit}(i_{\nu_n,\nu})\).
			
			\begin{claim}
				\(\vec{\kappa}\) generates the measure \(i_{\bar{\xi},\nu}(\bar{U})\).
				
				\begin{proof}
					Work in \(N_\nu\). Then we have
					\(A\in i_{\bar{\xi},\nu}(\bar{U})\) if and only if there is some \(n<\omega\) and \(\bar{A}\in N_{\nu_n}\) such that \(i_{\nu_n,\nu}(\bar{A})=A\) and \(\bar{A}\in i_{\bar{\xi},\nu_n}(\bar{U})\). Equivalently,  \(\kappa^{(n)}\in i_{\nu_n,\nu_{n+1}}(\bar{A})\) for some \(n<\omega\), i.e. \(\langle \kappa_n,\kappa_{n+1},\kappa_{n+2},...\rangle\in A\) for some \(n<\omega\).\\
					
					Note that \(K^M|(\kappa_\nu^{+})^{K^M}=N_\nu|(\kappa_\nu^{+})^{N_\nu}\), so \(i_{\bar{\xi},\nu}(\bar{U})\) is still a \(K^M\)-ultrafilter.
				\end{proof}
			\end{claim}
			
			So, since \(^{\omega}M\subseteq M\), we have that \(i_{\bar{\xi},\nu}(\bar{U})\in M\). Since \(\mathrm{crit}(i_{\bar{\xi},\nu})>\omega_2\), Lemma \ref{beginninglemma}(ii) implies that \(i_{\bar{\xi},\nu}(\bar{U})\in K^M=N_{\theta}\). But then \(i_{\bar{\xi},\nu}(\bar{U})\in N_{\nu+1}\), which is a contradiction. Hence, there is an inner model with a measurable cardinal \(\lambda\) such that \(o(\lambda)=\lambda^{++}\).
		\end{proof}
	\end{theorem}

\end{document}